\documentclass[11pt]{article}
\usepackage{graphicx}
\usepackage{color}
\usepackage{amsmath,amsfonts,amssymb,amsthm}
\def\ds{\displaystyle}
\def\forall{\hbox{for all}~}
\def\L{{\bf L}}

\def\Nat{{\mathbb N}}
\def\ve{\varepsilon}
\def\n{\noindent}

\def\R{{\mathbb R}}

\def\vp{\varphi}

\def\tv{\mathop{\mathrm{Tot.Var.}}}
\def\vs{\vskip 2em}

\def\v{\vskip 1em}
\def\D{{\cal D}}
\def\O{{\cal O}}

\def\begi{\begin{itemize}}
\def\endi{\end{itemize}}
\def\C{{\cal C}}
\def\dint{\int\!\!\int}
\def\ov{\overline}
\def\Tilde{\widetilde}
\def\Hat{\widehat}
\def\bega{\begin{array}}
\def\enda{\end{array}}
\def\meas{\hbox{meas}}
\def\bel{\begin{equation}\label}
\def\eeq{\end{equation}}
\def\sqr#1#2{\vbox{\hrule height .#2pt
\hbox{\vrule width .#2pt height #1pt \kern #1pt
\vrule width .#2pt}\hrule height .#2pt }}
\def\square{\sqr74}
\def\endproof{\hphantom{MM}\hfill\llap{$\square$}\goodbreak}

\newtheorem{theorem}{Theorem}[section]
\newtheorem{corollary}{Corollary}[section]
\newtheorem{lemma}{Lemma}[section]

\begin{document}
\title{\bf  Unique Solutions to Hyperbolic Conservation Laws with a Strictly Convex Entropy}
\vs

\author{Alberto Bressan$^*$ and  Graziano Guerra$^{**}$
\\
\, \\
{\small (*) Department of Mathematics, Penn State University,}\,\\ 
\small (**) Department of Mathematics and Applications,
  University of Milano - Bicocca. \\ \, \\
{\small E-mails:  axb62@psu.edu, ~ graziano.guerra@unimib.it}
\,\\}
\maketitle

\begin{abstract} Consider a strictly hyperbolic $n\times n$ system of conservation laws,
where each characteristic field is either genuinely nonlinear or linearly degenerate.
In this standard setting, it is well known that there exists a Lipschitz semigroup of 
weak solutions, defined on a domain of functions with small total variation.
If the system admits a strictly convex entropy, we give a short proof 
that every entropy weak solution taking values within the domain of the semigroup
coincides with a semigroup trajectory.  
The result shows that the assumptions of 
``Tame Variation" or ``Tame Oscillation", previously used to achieve uniqueness,
can be removed in the presence of a strictly convex entropy.

\end{abstract}
 \v
\section{Introduction}
\label{sec:1}
\setcounter{equation}{0}
We consider the Cauchy problem for  a strictly hyperbolic $n\times n$ system of conservation laws in one space dimension:
\bel{1} u_t + f(u)_x ~=~0,\eeq
\bel{2} u(0,x)~=~\bar u(x).\eeq
As usual, $f:\Omega\to\mathbb{R}^{n}$ is the flux, defined on some open set  $\Omega\subset \mathbb{R}^{n}$. We assume that each characteristic family is either genuinely nonlinear or linearly degenerate.
In this setting, it is well known that there exists a Lipschitz continuous semigroup $S:\D\times [0, +\infty[\, \mapsto \D$ of entropy weak solutions, defined on a domain $\D\subset\L^1(\R;\,\R^n)$ containing all functions
with sufficiently small total variation \cite{Bbook, BC95, BCP,  BLY, D2, HR}. The trajectories of this 
semigroup are  the unique limits of  front tracking approximations, and also of Glimm  
approximations \cite{BGlimm} and of
vanishing viscosity approximations \cite{BiB}.
We recall that the semigroup is globally Lipschitz continuous w.r.t.~the $\L^1$ distance. Namely,
there exists a constant $L$ such that
\bel{Slip}
\bigl\| S_t\bar u - S_s\bar v\bigr\|_{\L^1}~\leq~L\Big( |t-s| + \|\bar u - \bar v\|_{\L^1}\Big)\qquad
\forall s,t\geq 0,~~\bar u,\bar v\in\D.\eeq

Given any weak solution $u=u(t,x)$ of (\ref{1})-(\ref{2}),  
various conditions have been derived in \cite{BG, BLF, BL} which guarantee
the identity  
\bel{uniq} u(t)~=~S_t\bar u\qquad\forall t\geq 0.\eeq
Since the semigroup $S$ is unique, the identity
(\ref{uniq}) yields  the uniqueness of  solutions to the Cauchy problem (\ref{1})-(\ref{2}).
In addition to the standard assumptions, earlier results required some additional
regularity conditions, such as ``Tame Variation" or ``Tame Oscillation", controlling
the behavior of the solution near a point where the variation is locally small.

Aim of the present note is to show that, if the $n\times n$ system (\ref{1}) is endowed with a strictly convex entropy $\eta(\cdot)$, then every entropy-weak solution $t\mapsto u(t)$ 
taking values within the domain $\D$
of the semigroup satisfies (\ref{uniq}).  In other words, uniqueness is guaranteed without
any further regularity assumption.

As in \cite{BG, BLF, BL}, 
the proof relies on the elementary error estimate
\bel{errest}
\bigl\| u(t)- S_t\bar u\bigr\|_{\L^1}~\leq~L\cdot\int_0^t \liminf_{h\to 0+}\,
{\bigl\|u(\tau+h)- S_h u(\tau)\bigr\|_{\L^1}\over h}\, d\tau\,.\eeq
Assuming that the system is endowed with a strictly convex entropy, 
we will prove that the integrand is zero for a.e.~time $\tau\geq 0$. 
Following an argument introduced in \cite{BGlimm}, this is achieved by two estimates:
\begi
\item[(i)] In a neighborhood of a point $(\tau,y)$ where $u(\tau,\cdot)$ has a large jump,
the weak solution $u$ is compared with the solution to a Riemann problem.
\item[(ii)] In a region where the total variation is small,
the weak solution $u$ is compared with the solution to a linear system with constant
coefficients.\endi

The main difference is that here we estimate the lim-inf in (\ref{errest}) only at times
$\tau$ which are Lebesgue points for a countable family of total variation functions 
$W^{\xi,\zeta}(\cdot)$, defined at  (\ref{Wpm}).

To precisely state the result, we begin by collecting the main assumptions.
\v\begi
\item[{\bf (A1)}] {\bf (Conservation equations)}  {\it The function $u=u(t,x)$
is a weak solution of the Cauchy problem (\ref{1})-(\ref{2}) taking values within the domain of the semigroup.\\
More precisely, $u:[0,T]\mapsto\D$ is continuous w.r.t.~the $\L^1$ distance.
The identity $u(0,\cdot)=\bar u$ holds  in $\L^1$, and moreover
\bel{4}\dint \big(u\varphi_t+f(u)\varphi_x\big) ~dxdt=0\eeq
for every $\C^1$ function $\varphi$ with compact support contained
inside the open strip $\, ]0,T[\,\times\R$.}
\endi
Regarding the entropy conditions, 
we assume that the system (\ref{1})
admits 
a $\C^2$ entropy function $\eta:\Omega\mapsto\R$ with entropy flux $q$,
so that the equality $\nabla q(\omega) = \nabla \eta (\omega) Df(\omega)$ holds for all $\omega \in
\Omega$. We also assume that the entropy $\eta$ satisfies the strict convexity condition
\bel{scvx}
\eta(\omega) ~\geq~\eta(\ov\omega) +\nabla\eta(\ov \omega) \cdot (\omega-\ov \omega) + c_0 |\omega-\ov\omega|^2,\eeq
for some $c_0>0$ and  every couple of states $\omega,\; \ov \omega\in\Omega$.
As usual, we say that a weak solution $u$ is entropy-admissible if it satisfies:
\begi
\item[{\bf (A2)}] {\bf (Entropy admissibility condition)} 
{\it 
For every $\C^1$ function $\varphi\geq 0$ with compact support contained
inside the open strip $\, ]0,T[\,\times\R$, one has}
\bel{8}\dint \big(\eta(u)\varphi_t+q(u)\varphi_x\big) ~dxdt~\geq~0.\eeq
\endi
Our result can be simply stated as:
\begin{theorem}\label{t:1} Let (\ref{1}) be a strictly hyperbolic $n\times n$ system, where
each characteristic field is either genuinely nonlinear or linearly degenerate, and which
admits a strictly convex entropy $\eta(\cdot)$  as in (\ref{scvx}). 
Then every entropy-weak solution $u:[0,T]\mapsto\D$,  taking values within the 
domain of the semigroup, coincides with a semigroup trajectory.  \end{theorem}

The theorem will be proved in Section~\ref{sec:3}.    We remark that, restricted to a class of 
$2\times 2$ systems, a more elaborate
proof of this result was recently given in \cite{CKV}.

In our view, the main interest in the above uniqueness theorem is that it opens the door to the 
study of uniform convergence rates for a very wide class of approximation algorithms.
This will be better explained in the concluding remarks, in Section~\ref{sec:4}.

\v

\section{Preliminary lemmas}  
\label{sec:2}
\setcounter{equation}{0}

Let $M$ be an upper bound on the total variation of all functions in the domain
$\mathcal{D}$ of the semigroup:
\begin{equation}
  \label{eq:boundtv}
  \tv\left\{u;\mathbb{R}\right\}\le M,\quad \text{ for all }u\in\mathcal{D}.
\end{equation}
Since by assumption our solution $u(t,\cdot)\in\mathcal{D}$,
for sake of definiteness we shall assume that it is right continuous, namely
$u(t,x)= \lim_{y\to x+} u(t, y)$.
By~\cite[Theorem~4.3.1]{D2}, we have
the Lipschitz bound
 \bel{ulip}
 \bigl\| u(t_2,\cdot)-u(t_1,\cdot)\bigr\|_{\L^1(\R)}~\leq~C_M\,(t_2-t_1)
 \qquad\qquad \forall 0\leq t_1\leq t_2\,,\eeq
 for some constant $C_M>0$ depending only on  $M$ and on the flux
 $f$.

We begin by reviewing the well known fact  that the 
entropy has finite propagation speed.
\begin{lemma}
  \label{l:21}
  Let $u=u(t,x)$ be a function satisfying {\bf (A1)} and {\bf (A2)}.
  Then there exists two constants $\Hat C,\;\hat \lambda >0$ such
  that the following holds.   For any constant state $u^{*}\in \Omega$,  any
  $a<b$, and any $0\le\tau<\tau'$ with
  $2\hat \lambda (\tau' - \tau )< b-a$, 
one has
  \begin{equation}
    \label{eq:loc}
    \int_{a+\hat \lambda \left(\tau'-\tau\right)}^{b-\hat \lambda \left(\tau'-\tau\right)}\left|u(\tau',x) - u^{*}\right|^{2}\; dx
   ~ \le~ \Hat C \int_{a}^{b}\left|u(\tau,x) - u^{*}\right|^{2}\; dx.
  \end{equation}
\end{lemma}

{\bf Proof.}
Given the constant state $u^*\in\Omega$, for all $\omega\in\Omega$ 
define the relative entropy
  $\eta\left(\omega \mid u^{*}\right)$ and the corresponding entropy flux
  $q\left(\omega \mid u^{*}\right)$ as
 \bel{etq}\bega{rl}
    \eta\left(\omega \mid u^{*}\right) &=~ \eta\left(\omega \right) -
    \eta\left(u^{*}\right)-\nabla \eta\left(u^{*}\right)\left(\omega -u^{*}\right),\\[2mm]
    q\left(\omega \mid u^{*}\right) &=~ q\left(\omega \right) -
    q\left(u^{*}\right)-\nabla \eta\left(u^{*}\right)\bigl(f\left(\omega \right)-f\left(u^{*}\right)\bigr).
  \enda\eeq
  The equations \eqref{4} and \eqref{8} yield
  \begin{equation}
    \label{eq:relEnt}
      \eta\left(u \mid u^{*}\right)_{t}+  q\left(u \mid u^{*}\right)_{x}~\le~ 0,
    \end{equation}
    while~\eqref{scvx} implies
    \begin{equation}
      \label{eq:etaest}
      \eta\left(\omega \mid u^{*}\right)~\ge~ c_{0}\left|\omega - u^{*}\right|^{2},\qquad     \text{ for all }~ \omega, u^*\in\Omega.
    \end{equation}
By the $\C^2$ regularity of the functions $\eta,q$, there exists a constant $C'$
such that
    \begin{equation}
      \label{eq:Cest}
\eta\left(\omega \mid u^{*}\right)~\le~ C' \left|\omega
  - u^{*}\right|^{2},\qquad q\left(\omega \mid u^{*}\right)~\le~ C' \left|\omega
  - u^{*}\right|^{2},\quad
    \text{ for all }~ \omega, u^*\in\Omega.
\end{equation}
In view of (\ref{eq:etaest}),  there exists a sufficiently large constant
  $\hat \lambda >0$ such that
  \begin{equation}
    \label{eq:sidesest}
    -\hat \lambda \eta\left(\omega\mid u^{*}\right)\pm q\left(\omega\mid u^{*}\right)~\le~ 0,\qquad
    \text{ for all }~ \omega, u^*\in\Omega.
  \end{equation}
Using \eqref{eq:relEnt}, in connection with
test functions that approximate the
  characteristic function of the trapezoid $$\Gamma~=~\bigl\{(t,x) \,;~ \tau<t<\tau',
   ~~ a+\hat \lambda (t-\tau) < x < b-\hat \lambda (t-\tau)\bigr\},$$ and
 recalling \eqref{eq:sidesest}, we obtain
  \begin{displaymath}
    \begin{split}
      \int_{a+\hat \lambda (\tau'-\tau)}^{b-\hat \lambda (\tau'-\tau)}\eta\left(u\mid
        u^{*}\right)\left(\tau',x\right)\; dx&~\le ~
      \int_{a}^{b}\eta\left(u\mid
        u^{*}\right)\left(\tau,x\right)\; dx \\
      ~&\quad +\int_{\tau}^{\tau'}\left(-\hat \lambda \eta\left(u\mid u^{*}\right) + q\left(u\mid u^{*}\right)\right)\left(t,a+\hat
        \lambda \left(t-\tau\right)\right)\; dt\\
      &\quad +\int_{\tau}^{\tau'}\left(-\hat \lambda \eta\left(u\mid u^{*}\right) - q\left(u\mid u^{*}\right)\right)\left(t,b-\hat
        \lambda \left(t-\tau\right)\right)\; dt\\
      &~\le~ 
      \int_{a}^{b}\eta\left(u\mid
        u^{*}\right)\left(\tau,x\right)\; dx.
    \end{split}
  \end{displaymath}
  Together with~\eqref{eq:etaest}-\eqref{eq:Cest}, this proves the lemma.  
\endproof

Throughout the following, without loss of generality we shall always assume $\hat\lambda=1$.
We observe that this can always be achieved by a suitable rescaling of the time variable:
$$\tilde t~=~\kappa  t.$$
Similarly,  to simplify notation, we also assume that all wave speeds lie in the interval $[-1,1]$.

Given any $\tau\geq 0$ and any bounded interval $]a,b[$ with $-\infty\le a < b \le +\infty$, we  consider the open  intervals
\bel{Jt} J(t)~=~\bigl] a+ {}(t-\tau)\,,~b- {}(t-\tau)\bigr[\,,\qquad\qquad \tau\le t< \tau+
\frac{b-a}{2}\,.\eeq

Toward the proof of Theorem~\ref{t:1}, in order  to replace the ``Tame Variation" condition,
the main tool is provided by the following elementary lemma.
\begin{lemma} \label{l:22} 
In the setting of Theorem~\ref{t:1}, for some constant $C>0$ the following holds.
Let $u=u(t,x)$ be any entropy weak solution to (\ref{1}).  Then

\bel{ele}
\int_{J(t)} \bigl| u(t,x) - u(\tau, x)\bigr|\, dx~\leq~C (t-\tau) \cdot\tv\bigl\{ u(\tau,\cdot)\,;~]a,b[\bigr\}.
\eeq
\end{lemma}

\begin{figure}[htbp]
\centering
 \includegraphics[scale=0.5]{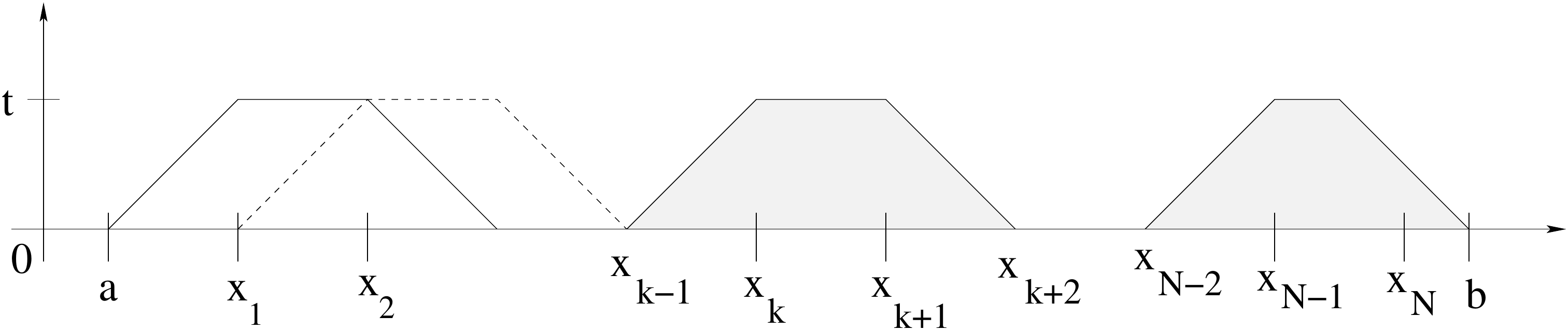}
    \caption{\small  The covering of the interval $[a,b]$ used in the proof of Lemma~\ref{l:22}.}
\label{f:hyp247}
\end{figure}

{\bf Proof.}  {\bf 1.}
 We first consider the case where $-\infty <a < b<+\infty$. 
 For notational simplicity, w.l.o.g.~we assume that $\tau=0$.  Given a time
  $0<t<\frac{b-a}{2 }$, as shown in Fig.~\ref{f:hyp247}
  we define the points $x_k$, the values $u_k$ and the integer $N\geq 1$ such that
  \bel{xkdef}
  x_k~=~ a+kt,\qquad\qquad  u_k = u(0, x_k),\qquad\qquad x_N\leq b< x_{N+1}\,.\eeq
For $k=1,2,\ldots, N-2$, we apply Lemma~\ref{l:21} with $u^*=u_k$ on the trapezoids
$$\Gamma_k~\doteq~\bigl\{ (s,x)\,;~~s\in [0,t]\,,~~x_{k-1}+s<x<x_{k+2}-s\bigr\}.$$
Then we apply the same lemma with $u^*= u_{N-1}$ on the domain
$$\Tilde \Gamma_N~\doteq~\bigl\{ (s,x)\,;~~s\in [0,t]\,,~~x_{N-2}+t<x<b-t\bigr\}.$$

This yields the estimates
\bel{b11}\bega{rl} 
\ds      \int_{x_{k}}^{x_{k+1}} \bigl| u(t,x)- u_{k}\bigr|^2\, dx
      ~&\leq\ds~\Hat C\,\int_{x_{k-1}}^{x_{k+2}} \bigl| u(0,x)- u_{k}\bigr|^2\, dx\\[4mm]
      &\leq~\Hat C \cdot 3t \cdot \Big( \tv\{
      u(0,\cdot)\,;~]x_{k-1},
     x_{k+2}[\bigr\} \Big)^2,
    \enda
 \eeq
 \bel{b22}  \int_{x_{N-1}}^{b-t} \bigl| u(t,x)- u_{N-1}\bigr|^2\, dx
      ~\leq\ds~\Hat C \cdot 3t \cdot \Big( \tv\{
      u(0,\cdot)\,;~]x_{N-2},
     b[\bigr\} \Big)^2.\eeq
\v
{\bf 2.} Define the piecewise constant approximation $\bar u:[a +   t,b-  t[\,\mapsto 
  \Omega$ by setting
  \begin{equation}
    \label{eq:defbaru}
    \bar u(x)\,=\,
    u_{k}\qquad \text{ if }~x \in \bigl[x_{k},x_{k+1}\bigr[\,,\qquad\qquad  k = 1,\ldots, N-1.
  \end{equation}
  Using Cauchy's inequality and the bounds (\ref{b11})-(\ref{b22}),  we obtain
  \bel{est3}
    \bega{l}\ds
      \int_{a+  t}^{b-  t} \bigl|u(t,x)-\bar u(x)\bigl|\,
      dx~ =~\sum_{k=1}^{N-2}\int_{x_k}^{x_{k+1}} \bigl| u(t,x)-
      u_{k}\bigr|\, dx+ \int_{x_{N-1}}^{b-t} \bigl| u(t,x)-
      u_{N-1}\bigr|\, dx\\[4mm]
      \qquad \ds
      \le~ \sum_{k=1}^{N-2}\sqrt{t}\left(\int_{x_{k}}^{x_{k+1}} \bigl| u(t,x)- u_{k-1}\bigr|^{2}\, dx\right)^{1/2}+\sqrt t  \left(\int_{x_{N-1}}^{b-t} \bigl| u(t,x)-
      u_{N-1}\bigr|^2\, dx\right)^{1/2}  \\[4mm]
   \qquad \ds   \leq ~ \sqrt t\cdot \sqrt{\Hat C\cdot 3t }\cdot
  \left(    \sum_{k=1}^{N-2}\tv\{ u(0,\cdot)\,;~]x_{k-1}, x_{k+2}[\bigr\} +
  \tv\{
      u(0,\cdot)\,;~]x_{N-2},
     b[\bigr\} \right).
    \enda
    \eeq
    Observing that every point $x\in [a,b]$ is contained in at most three 
    open intervals 
    $]x_{k-1}, x_{k+2}[\,$, from (\ref{est3}) we conclude
   \bel{est4}
    \int_{a+  t}^{b-  t} \bigl|u(t,x)-\bar u(x)\bigl|\, dx~
    \leq ~ \sqrt{3\Hat C}\cdot 3   t\cdot \tv\bigl\{ u(0,\cdot)\,;~]a,
    b[\bigr\}.
  \eeq
  \v
  {\bf 3.} Next, we compute
  \bel{est5}
  \bega{l}\ds
      \int_{a+  t}^{b-  t} \bigl|u(0,x)-\bar u(x)\bigl|\,
      dx~=\ds~ \sum_{k=1}^{N-2}\int_{x_k}^{x_{k+1}} \bigl| u(0,x)-
      u_{k}\bigr|\, dx +    \int_{x_{N-1}}^{b-t}  \bigl| u(0,x)-u_{N-1}\bigr|\, dx\\[4mm]
      \quad \leq \ds~ \sum_{k=1}^{N-2} t\cdot \tv\bigl\{ u(0,\cdot)\,;~
      ]x_k, x_{k+1}[\,\bigr\} +    t\cdot \tv\bigl\{ u(0,\cdot)\,;~
      ]x_{N-1}, b-t[\,\bigr\}  \\[4mm]
      \quad \ds \le ~  t \cdot \tv\bigl\{ u(0,\cdot)\,;~]a, b[\bigr\}.
    \enda
  \eeq
  Combining (\ref{est4}) with (\ref{est5}) we obtain a proof of the
  lemma for finite $a$ and $b$.
  
  Letting $a\to -\infty$ or $b\to +\infty$ we see that 
  the same conclusion  remains valid also for unbounded intervals,
  such as $]-\infty, b[\,$ or $\,]a, +\infty[\,$.
  \endproof

\section{Proof of the theorem}  
\label{sec:3}
\setcounter{equation}{0}

 We are now ready to give a 
 proof of Theorem~\ref{t:1}, in several steps.
 \v
 {\bf 1.} By the structure theorem for BV functions \cite{AFP, EG}, there is a null set of times
${\cal N}\subset [0,T]$ such that the following holds.

Every point $(\tau,\xi)\in [0,T]\times \R$ with $\tau\notin {\cal N}$ is either a point of approximate continuity, or a point of approximate jump of the function $u$.  
In this second case, there exists states 
$u^-,u^+\in\R^n$ and a speed $\lambda\in\R$ such that, calling
\bel{5}U(t,x)~\doteq~
\left\{\bega{rl} u^-\qquad &\hbox{if}\qquad \left(x-\xi\right)<\lambda \left(t-\tau\right),
\\[2mm] u^+\qquad &\hbox{if}\qquad \left(x-\xi\right)>\lambda \left(t-\tau\right),\enda\right.
\eeq
there holds
\bel{6}\lim_{r\to 0+}~{1\over r^2}
\int_{-r}^r\int_{-r}^r\Big| u(\tau+t,~\xi+x)-U(\tau + t,\xi + x)\Big|~dxdt~=~0.\eeq
The conservation equations (\ref{4}) imply that  
the piecewise constant function $U$ must be a
weak solution to the system of conservation laws (see~\cite[Theorem~4.1]{Bbook}),
satisfying the Rankine-Hugoniot equations:
\bel{RH}
f(u^+) - f(u^-)~=~\lambda (u^+-u^-).\eeq
Moreover, the entropy condition (\ref{8}) implies
\bel{qeta}
q(u^+) - q(u^-) ~\leq~\lambda\bigl(\eta(u^+) - \eta(u^-)\bigr).\eeq

Next, we observe that, for every couple of rational points $\xi,\zeta\in {\mathbb Q}$, the 
scalar function
\bel{Wpm} W^{\xi,\zeta} (t)~\doteq~\left\{ \bega{cl} \tv\big\{ u(t)\,;~]\xi+t\,,~\zeta-t[\,\bigr\}
\qquad &\hbox{if}\quad \xi+t<\zeta-t\,,\\[2mm]
0\qquad &\hbox{otherwise,}\enda\right.\eeq
is bounded and measurable (indeed, it is lower semicontinuous).  Therefore a.e.~$t\in [0,T]$ is a Lebesgue point.
We denote by ${\cal N}'\subset [0,T]$ the set of all times $t$ which are NOT
Lebesgue for at least one of the countably many functions  $W^{\xi,\zeta}$.
Of course, ${\cal N}'$ has zero Lebesgue measure.

In view of (\ref{errest}), we will prove the theorem by establishing  the following claim.
\begi
\item[{\bf (C)}] {\it For every $\tau\in [0,T]\setminus ({\cal N\cup \cal N'})$ and $\ve>0$,
one has}
\bel{limsup}
\limsup_{h\to 0+} {1\over h} \Big\| u( \tau + h) - S_h u( \tau)\Big\|_{\L^1}~\leq~\ve.\eeq
\endi 
 \v
 {\bf 2.} Assume $\tau \notin \cal N\cup \cal N'$.  
 Since $u(\tau,\cdot)$ has bounded variation, 
we define points
\begin{displaymath}
-\infty = y_{-1}<y_0<y_1<\cdots<y_N<y_{N+1}=+\infty, 
\end{displaymath}
\begin{displaymath}
y_{k+1} ~=~\sup~\Big\{ x>y_{k}\,;~~\tv\bigl\{ u(\tau,\cdot)\,;~~]y_k, x[\,\bigr\} \leq\ve\Big\}.
\end{displaymath}
Since $u(\tau,\cdot)$ is right continuous, we have
\begin{equation}
  \label{tvep}
  \begin{cases}
    \tv\bigl\{ u(\tau, \cdot)\,;~]y_{k-1}, y_k[\,\bigr\}~\leq~\ve,&\text{ for }k=0,\ldots,N+1,\\
    \tv\bigl\{
    u(\tau, \cdot)\,;~]y_{k-1}, y_k]\,\bigr\}~\ge~\ve,&\text{ for }k=0,\ldots,N,\\
    N~\leq~{M\over \ve}.
  \end{cases}
\end{equation}
Where $M$ is an upper bound for the total variation of all functions $\bar u\in\D$, as in 
(\ref{eq:boundtv}).

Then we choose points $y_k', y_k''$ such that
$$-\infty<y_0<y_0''~\leq~y_1'<y_1<y_1''~\leq~ y_2'<y_2<y_2''~<~\cdots~\leq~y_N'<y_N<+\infty$$
\begin{equation}
  \label{ykep}
  \left\{
\bega{l}\tv\bigl\{u( \tau,\cdot) \,;~]y'_k, y_k[\,\bigr\}\,\leq \,\ve^2,\quad k=1,\ldots N,\\[3mm]
\tv\bigl\{u(\tau,\cdot) \,;~]y_k, y''_k[\,\bigr\}\leq
\,\ve^2,\,\quad k=0,\ldots N-1,\\[3mm]
\text{ all values $y_k'-\tau$, $y''_k-\tau$ are rational}.\enda\right.
\end{equation}

\v
\begin{figure}[htbp]
\centering
 \includegraphics[scale=0.5]{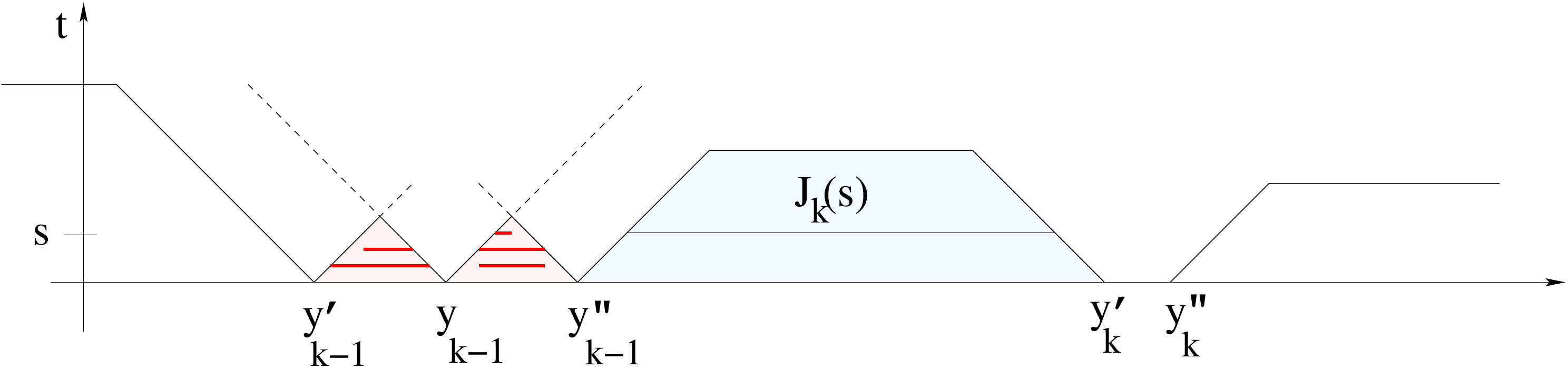}
    \caption{\small  The points $y'_k<y_k<y''_k$ constructed 
 in the proof of the theorem. Typically, $y_k$ is the location of a shock. 
Since $y_k\pm \tau$ need not be rational, the additional points $y_k', y_k''$ must be considered.}
\label{f:hyp251}
\end{figure}

\n{\bf 3.}  For any given $y\in \R$, we denote by $U^\sharp=
U^\sharp_{(u,\tau,y)}(t,x)$ the solution, for  $t\ge \tau$,  to the 
Riemann problem for (\ref{1}) with initial data at $t=\tau$:
\bel{Rdata}\bar u( x)~=~\left\{ \bega{rl} u(\tau, y-)\qquad &\hbox{if} ~~x< y,\cr
u(\tau, y+)\qquad &\hbox{if} ~~x> y.\enda\right.\eeq
Moreover, for every given $k=1,\ldots,N$ we denote by $U^\flat= U^\flat_{(u,\tau,k)}(t,x)$ the solution to the 
linear Cauchy problem with constant coefficients
\bel{linh} v_t + A v_x~=~0,\qquad \qquad  v(\tau,x) = u(\tau, x).\eeq
Here the $n\times n$ matrix $A$ is the Jacobian matrix of $f$ computed
at the midpoint of the interval $[y''_{k-1}, y'_k]$. Namely, 
$$~A~=~ Df\left(\Tilde u_{k}\right),
\qquad \Tilde u_{k}=u\left(\tau, {y''_{k-1}+ y'_k\over 2}\right),\quad k=1,\ldots,N.$$
 With reference to Fig.~\ref{f:hyp251}, to estimate  the lim-sup
in (\ref{limsup}), we need to estimate three types of integrals.

\begi
\item[(I)] The integral of $\bigl| u(t,x) - U^\sharp _{\left(u;\tau,y\right)}(t,x)\bigr|$ over the interval
$$ J_y(t)~\doteq~
\bigl[  y -(t-\tau)\,,~ y +(t-\tau) \bigr],$$
for all points $y\in \bigl\{ y_0, y_0'', y_1', y_1, y_1'',\ldots, y_N', y_N\bigr\}$.

\item[(II)]The integral of $\bigl| u(t,x) - U^\flat_{\left(u;\tau,k\right)} (t,x)\bigr|$  over the interval
\bel{Jk}J_{k}(t)~=~\bigl] y_{k-1}''+ (t-\tau)\,,~
 y_{k}'- (t -\tau) \bigr[,\quad k=1,\ldots,N.\eeq

\item[(III)] The integral of $\bigl| u(t,x) - u(\tau,x)\bigr|$ over the intervals
$$\left\{
\bega{rl}J_k'(t)&=~\bigl]  y_k' +(t-\tau) \,,~ y_k -(t-\tau) \bigr[,\quad k=1,\ldots,N,\\[2mm]
J''_k(t)&=~\bigl]  y_k +(t-\tau) \,,~ y''_k -(t-\tau) \bigr[,\quad k=0,\ldots,N-1,\\[2mm]
J_{0}(t)&\dot=~\left]-\infty,y_{0}-\left(t-\tau\right)\right[,\\[2mm]
J_{N+1}&\dot=~\left]y_{N}+\left(t-\tau\right),+\infty\right[.
\enda\right.
$$

\endi
\v
{\bf 4.}  To estimate integrals of type (I), assuming that $(\tau, y)$ is either
a Lebesgue point or  a point of approximate jump of the 
function $u$, we obtain
\bel{Ush}\lim_{h\to 0+} {1\over h}\int_{y-h}
^{y+h}\Big| u(\tau+h,~x)- U^\sharp
_{(u;\tau,y)}(\tau+h, x)\Big|~dx~=~0.\eeq

Indeed, by~\cite[Theorem~2.6]{Bbook}, setting
$u^{\pm}=u\left(\tau,y\pm\right)$, the function $U$
defined in~\eqref{5}  satisfies~\eqref{6} and consequently $u^{\pm}$
satisfy~\eqref{RH} and~\eqref{qeta}. It implies that
$u\left(\tau,y\pm\right)$ are connected by a single entropic shock whose
speed is $\lambda$. Consequently $U^{\sharp}_{\left(u;\tau,y\right)}$ coincides with
the piecewise constant 
function $U$ defined
in~\eqref{5} so that~\eqref{Ush} follows
from~\cite[Theorem~2.6]{Bbook}.

\v
{\bf 5.} We now estimate the integrals of type (II).
By construction, both values $y_{k-1}''-\tau$ and $y'_k + \tau$ are rational. Hence
$$y_{k-1}''\,=\,\xi+\tau,\qquad\qquad y_k'\,=\,\zeta-\tau,$$
for some $\xi, \zeta\in {\mathbb Q}$.
This implies that $\tau$ is a Lebesgue point of the map 
\bel{Vt} t~\mapsto~ V(t)~=~W^{\xi,\zeta}(t).\eeq

Let $\tilde\lambda_i=\lambda_{i}\left(\tilde u_{k}\right)$, $\tilde
l_i=l_{i}\left(\tilde u_{k}\right)$,
$\tilde r_i=r_{i}\left(\tilde u_{k}\right)$, $i=1,\ldots,n$, be
respectively the $i$-th eigenvalues and left and right eigenvectors of
the matrix $ A\doteq Df\big(\Tilde u_{k}\big)$.
We thus have
$$\tilde l_i\cdot U^\flat(t,x)~=~\tilde l_i\cdot U^\flat\big(\tau,
x-(t-\tau)\tilde\lambda_i\big)~=~\tilde l_i\cdot u\big(\tau,
x-(t-\tau)\tilde\lambda_i\big).$$
Following the proof of~\cite[Theorem~9.4]{Bbook}, fix any two points
$\zeta',\zeta''\in J_{k}(t)$, $\zeta'<\zeta''$ and consider the quantity
\bel{53}\bega{rl} E_i(\zeta',\zeta'')&\ds
\doteq ~ \tilde l_i\cdot\int_{\zeta'}^{\zeta''}\big(u(t,x)
-U^\flat(t,x)\big)~dx\\[3mm]
&=~ \ds\tilde l_i\cdot\int_{\zeta'}^{\zeta''}\Big(u(t,x)
-u\big(\tau,~x-(t-\tau)\tilde\lambda_i\big)\Big)~dx
.\enda \eeq

\begin{figure}[htbp]
\centering
 \includegraphics[scale=0.4]{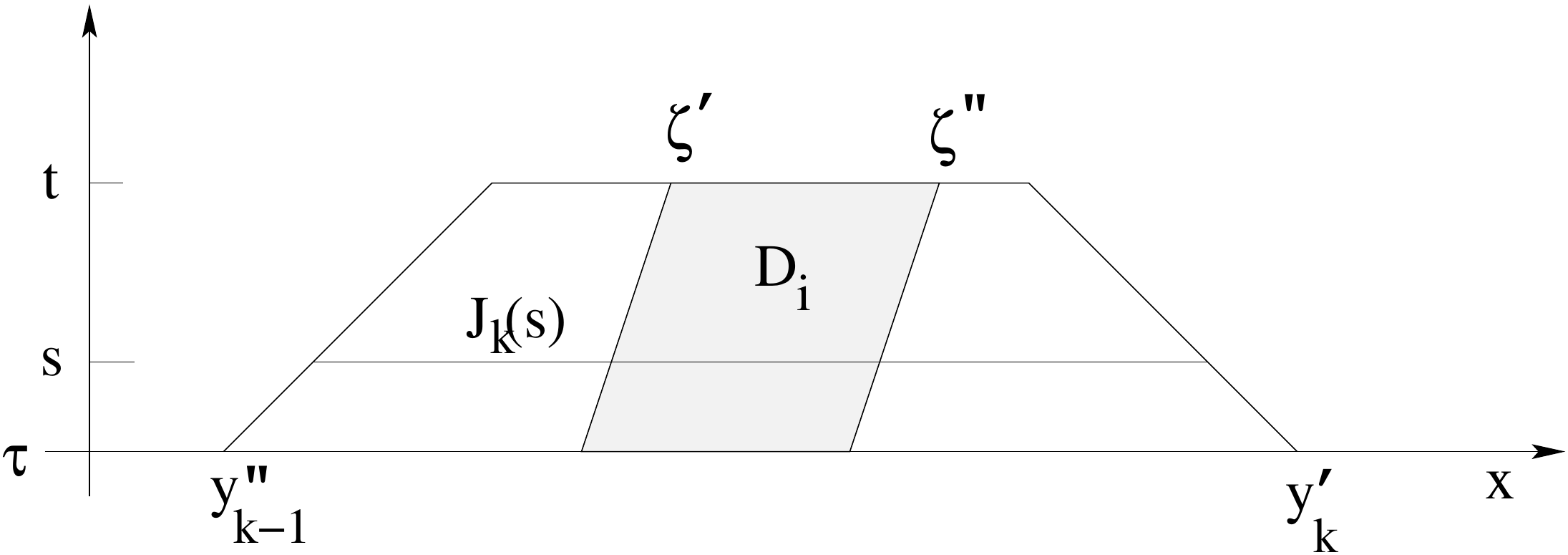}
    \caption{\small  The domain $D_i$ considered at (\ref{Di}).}
\label{f:hyp254}
\end{figure}

We apply the divergence theorem to the vector
$\big(u,\,f(u)\big)$
on the domain
\bel{Di}D_i~\doteq~\Big\{ (s,x);~~s\in[\tau,~t],~~
\zeta'-(t-s)\tilde\lambda_i\leq x\leq \zeta''-(t-s)\tilde
\lambda_i\Big\},\eeq
shown in Fig.~\ref{f:hyp254}.
Since $u$ satisfies the conservation equation (\ref{1}),
the difference between the integral of $u$ at the top and at the
bottom of the domain $D_i$ is thus
measured by the inflow from the left side minus
the outflow from the right side of $D_i$.
From (\ref{53}) it thus follows
\bel{54}\bega{rl} E_i(\zeta',\zeta'')=& ~\ds
\int_\tau^{t} \tilde l_i\cdot\Big(\big(f(u)-\tilde\lambda_i 
u\big)(s,~\zeta'-(t-s)\tilde \lambda_i)\Big) \,ds\cr
&\ds \qquad-\int_\tau^{t} \tilde l_i\cdot\Big(\big(f(u)-\tilde\lambda_i 
u\big)(s,~\zeta''-(t-s)\tilde \lambda_i)\Big)\,ds\cr
=&\ds \int_\tau^{t}  l_i(\tilde u_k)\cdot\Big(\big(f(u'(s))-\tilde
\lambda_i\left(\tilde u_k\right) 
u'(s)\big)-  \big(f(u''(s))-\lambda_i(\tilde u_k) 
u''(s)\big)
\Big)\,ds\cr
=&\ds \int_\tau^{t}  H\left(\Tilde u_k,u'(s),u''(s)\right)\,ds,
\enda \eeq
where we set
$$u'(s)\doteq u\big(s,~\zeta'-(t-s)\tilde\lambda_i\big),\qquad
u''(s)\doteq u\big(s,~\zeta''-(t-s)\tilde\lambda_i\big),$$
\begin{displaymath}
  H\left( u,u_{1},u_{2}\right)
  \,\doteq\,l_i( u)\cdot\Big(\big(f(u_{1})-\lambda_i(u) 
u_{1}\big)-  \big(f(u_{2})- \lambda_i (u)
u_{2}\big)
\Big).
\end{displaymath}
Observing that
\begin{itemize}
\item
  $H\left(u,u_{2},u_{2}\right)=0$,
\item
  $D_{u_{1}}H\left(u,u_{1},u_{2}\right)=l_{i}(u)\cdot
  \left(Df(u_{1})-\lambda_{i}(u)I\right)$,
\item
  $D_{u_{1}}H\left(u,u,u_{2}\right)=l_{i}(u)\cdot
  \left(Df(u)-\lambda_{i}(u)I\right)=0$,
\end{itemize}
we  estimate
$$
  \bega{rl}   H\left(u,u_{1},u_{2}\right)&=~
    H\left(u,u_{1},u_{2}\right)-
    H\left(u,u_{2},u_{2}\right)\\[4mm]
    &=~\ds
    \int_{0}^{1}D_{u_{1}}H\left(u,u_{2}+\sigma \left(u_{1}-u_{2}
      \right),u_{2}\right)\; d \sigma \cdot \left(u_{1}-u_{2}\right)\\[4mm]
    &=~\ds
    \int_{0}^{1}\left[D_{u_{1}}H\left(u,u_{2}+\sigma \left(u_{1}-u_{2}
      \right),u_{2}\right)-D_{u_{1}}H\left(u,u,u_{2}\right)\right]\; d
  \sigma \cdot \left(u_{1}-u_{2}\right)\\[4mm]
  &=~\ds \O(1)\cdot \bigl(\left|u_{1}-u\right|+\left|u_2-u\right|\bigr)\cdot\left|u_{1}-u_{2}\right|.
  \enda $$
Therefore,
\begin{displaymath}
    E_i(\zeta',\zeta'')~=~\O (1)\cdot \int_\tau^{t} \left|u'(s)-u''(s)\right|\cdot
    \Big(\left|u'(s)-\Tilde u_k\right|+ \left|u''(s)-\Tilde
        u_k\right|\Big)\, ds.
\end{displaymath}
Recalling (\ref{Jk}) and (\ref{Vt}), for any $x\in J_{k}(s)$ we now compute
\begin{displaymath}
    \left|u'(s)-\Tilde u_k\right| \le V(s)+
    \left|u\big(s,x\big)- u\left(\tau,x\right) \right| + \left| u\left(\tau,x\right)-\Tilde u_{k}\right|
    \le V(s)+
    \left|u\big(s,x\big)- u\left(\tau,x\right) \right| + V(\tau).
\end{displaymath}
Integrating w.r.t.~$x$ over the interval
$J_{k}(s)$, dividing by its length and using~\eqref{ulip}
 we obtain
\bel{gdef}\bega{rl}
    \left|u'(s)-\Tilde u_k\right| &\ds \le~ V(s)+V(\tau) +  \frac{1}{\meas \bigl(J_{k}(s)\bigr)}\int_{J_k(s)}\bigl|u(s,x)-u(\tau,x)\bigr|\, dx\\[4mm]
    &\ds =~V(s)+\ve +
    \O(1)\cdot
    \frac{C_M (s-\tau)}{\meas \bigl(J_{k}(s)\bigr)}~\dot=~g(s).
\enda \eeq
An entirely similar estimate clearly holds for $\left|u''(s)-\Tilde u_k\right|$. 
Hence
$$
  \bega{rl}
    E_i(\zeta',\zeta'')&\ds =~\O (1)\cdot \int_\tau^{t} \left|u'(s)-u''(s)\right|\cdot g(s)\,
    ds\\[4mm]
    &=~\ds \O(1)\cdot  \int_\tau^t \tv\Big\{ u(s)\,;~\bigl]\zeta'-(t-s)\tilde\lambda_{i},~
    \zeta''-(t-s)\tilde\lambda_{i}\bigr]\Big\}\cdot
    g(s)\, ds\\[4mm]
   &=~\O(1)\cdot\mu_{i}\bigl(]\zeta',\zeta'']\bigr).
\enda
$$
Here $\mu_i$ is the Borel measure defined by
$$
  \mu_{i}\bigl(]a,b[\bigr)~=~  \int_\tau^t \tv\Big\{ u(s)\,;~\bigl]a-(t-s)\tilde\lambda_{i},~
   b-(t-s)\tilde\lambda_{i}\bigr[\Big\}\cdot
    g(s)\, ds, $$   for any open interval $]a,b[ \,\subset J_{k}(t)$.  
    
According to Lemma~9.3 in \cite{Bbook}, we now have
\begin{displaymath}
  \begin{split}
    \int_{J_{k}(t)}\left|u(t,x)-U^{\flat}\left(t,x\right)\right|\; dx&=~ \O(1)\cdot
    \sum_{i=1}^{n} \int_{J_{k}(t)}\left|\tilde l_{i}\cdot
      \left(u(t,x)-U^{\flat}\left(t,x\right)\right)\right|\; dx\\
    &=~
    \O(1)\cdot \sum_{i=1}^{n}\mu_{i}\left(J_{k}(t)\right)~=~\O(1)\cdot
    \int_{\tau}^{t}V(s)\cdot g(s)\;
    ds.
  \end{split}
\end{displaymath}
In turn, this implies
\bel{uuf}
\frac{1}{t-\tau}  \int_{J_{k}(t)}\left|u(t,x)-U^{\flat}\left(t,x\right)\right|\;
dx~=~ \O(1)\cdot\left(\frac{\left\|g\right\|_{\infty}}{t-\tau}\int_{\tau}^{t}\left|V(s)-V(\tau)\right|\; ds
+\frac{V(\tau)}{t-\tau}\int_{\tau}^{t}g(s)\; ds\right).
\eeq

We now observe that, for all $s>\tau$ sufficiently close to $\tau$,
the function $g$ introduced at (\ref{gdef}) satisfies
\bel{gbound}
g(s)~\leq~V(s) +\ve + \O(1) \cdot \left(t-\tau\right).\eeq
Since  $t=\tau$ is a
Lebesgue point for $V$, taking the limit of (\ref{uuf}) as $t  \to\tau+$ we thus obtain
\bel{limuuf}
  \limsup_{t\to\tau+}
  \frac{1}{t-\tau}  \int_{J_{k}(t)}\left|u(t,x)-U^{\flat}\left(t,x\right)\right|\;
dx~=~ \O(1)\cdot V(\tau)\left(V(\tau)+\varepsilon\right)~=~ \O(1)\cdot\varepsilon^{2}.
\eeq

\v
{\bf 6.}  Finally, regarding integrals of type (III),  using Lemma~\ref{l:22} 
we obtain the bounds
\bel{J3}
\limsup_{h\to 0} ~{1\over h} \int_{y'_k+h}^{y_k-h}
\Big| u(\tau+h,~x)- u(\tau, x)\Big|~dx~=~\O(1)\cdot \ve^2,\eeq
\bel{J4}\limsup_{h\to 0} ~{1\over h} \int_{y_k+h}^{y''_k-h}
\Big| u(\tau+h,~x)- u(\tau, x)\Big|~dx~=~\O(1)\cdot \ve^2
,\eeq
\bel{J5}\limsup_{h\to 0} ~{1\over h} \int_{-\infty}^{y_0-h}
\Big| u(\tau+h,~x)- u(\tau, x)\Big|~dx~=~\O(1)\cdot \ve
,\eeq
\bel{J6}\limsup_{h\to 0} ~{1\over h} \int_{y_N+h}^{+\infty}
\Big| u(\tau+h,~x)- u(\tau, x)\Big|~dx~=~\O(1)\cdot \ve
.\eeq
\v
{\bf 7.} On the other hand, it is well known \cite{BGlimm, Bbook} that semigroup trajectories
satisfy entirely similar estimates. Indeed, at every point $y$ the difference
between the semigroup solution and the solution to a Riemann problem 
satisfies
\bel{S1}\lim_{h\to 0+} {1\over h}\int_{y-h}
^{y+h}\Big| \bigl(S_h u(\tau)\bigr) (x)- U^\sharp
_{(u;\tau,y)}(\tau+h, x)\Big|~dx~=~0.\eeq
Since the total variation of $u(\tau,\cdot)$ on the open interval $]y''_{k-1}, y_k'[$ is $\leq\ve$,
we have
\bel{S2}
\limsup_{h\to 0} ~{1\over h } \int_{y_{k-1}''+h}^{y_k'-h}
\Big|  \bigl(S_h u(\tau)\bigr) (x)- U^\flat_{(u,\tau,k)}(\tau+h, x)\Big|~dx~=~\O(1)\cdot \ve^2.\eeq
Moreover, since the total variation of $u(\tau,\cdot)$ on the open intervals 
$]y'_k, y_k[$ and $]y_k, y''_k[$ is $\leq\ve^2$,
we have
\bel{S3}
\limsup_{h\to 0} ~{1\over h} \int_{y'_k+h}^{y_k-h}
\Big|  \bigl(S_h u(\tau)\bigr) (x)- u(\tau, x)\Big|~dx~=~\O(1)\cdot \ve^2,\eeq
\bel{S4}\limsup_{h\to 0} ~{1\over h} \int_{y_k+h}^{y''_k-h}
\Big|  \bigl(S_h u(\tau)\bigr) (x)- u(\tau, x)\Big|~dx~=~\O(1)\cdot \ve^2,\eeq
and similarly
\bel{S5}
\limsup_{h\to 0} ~{1\over h} \int_{-\infty}^{y_0-h}
\Big|  \bigl(S_h u(\tau)\bigr) (x)- u(\tau, x)\Big|~dx~=~\O(1)\cdot \ve,\eeq
\bel{S6}
\limsup_{h\to 0} ~{1\over h} \int_{y_N+h}^{\infty}
\Big|  \bigl(S_h u(\tau)\bigr) (x)- u(\tau, x)\Big|~dx~=~\O(1)\cdot \ve.\eeq
\v
{\bf 8.}  Combining all the previous estimates, and recalling that the total number of 
intervals is $N\leq M\ve^{-1}$, we establish the limit (\ref{limsup}), proving the theorem.
\endproof

\section{Concluding remarks}  
\label{sec:4}
\setcounter{equation}{0}
The present analysis opens the door to the study of  
convergence and a posteriori error estimates
for a wide variety of approximation schemes.

Following \cite{BCS}, we say that $u=u(t,x)$ is an  $\ve$-approximate solution to (\ref{1}) if, 
given the time step $\ve = \Delta t$, the following holds.

\begi\item[{ \bf (AL)}] {\bf Approximate Lipschitz continuity.}  For every $\tau,\tau'\geq 0$
one has
$$\|u(\tau,\cdot)-u(\tau',\cdot)\|_{\L^1}~\leq~M\,\bigl(|\tau-\tau'|+\ve\bigr)\cdot\sup_{t\in 
[\tau, \tau']} \tv\bigl\{ u(t,\cdot)\bigr\}.
$$

\item[{\bf (P$_\ve$)}] {\bf Approximate conservation law, and approximate entropy inequality.}

{\it For every strip $[\tau,\tau' ]\times \R$ with $\tau,\tau'\in \ve\Nat$, and every test function $\vp\in \C^1_c(\R^2)$, there holds
\bel{wsol}\bega{l}\ds \left| \int u(\tau,x)\vp(\tau,x)\, dx-\int u(\tau',x)\vp(\tau',x)\, dx+
\int_\tau^{\tau'}\! \!\int \bigl\{ u\vp_t+f(u)\vp_x\bigr\}\,dx\,dt\right|\\[4mm]
\qquad\qquad\ds \leq~C \ve \|\vp\|_{W^{1,\infty}}
\cdot(\tau'-\tau)\cdot \sup_{t\in 
[\tau, \tau']} \tv\bigl\{ u(t,\cdot)\bigr\}.\enda\eeq

Moreover, given a uniformly convex entropy $\eta$ with flux $q$, assuming $\vp\geq 0$ 
one has the entropy inequality}
\bel{eiq}
\bega{l}\ds  \int \eta(u(\tau,x))\vp(\tau,x)\, dx-\int \eta(u(\tau',x))\vp(\tau',x)\, dx+
\int_\tau^{\tau'}\! \!\int \bigl\{\eta(u)\vp_t+q(u)\vp_x\bigr\}\,  dxdt
\\[4mm]\qquad\qquad\ds \geq~-C \ve \|\vp\|_{W^{1,\infty}}
\cdot(\tau'-\tau)\cdot \sup_{t\in 
[\tau, \tau']} \tv\bigl\{ u(t,\cdot)\bigr\}.\enda\eeq
\endi

In the above setting, the recent paper \cite{BCS} has established {\bf a posteriori}
error estimates,
assuming that the total variation of $u(t,\cdot)$ remains small, so that $u(t,\cdot)$ remains 
inside the domain of the semigroup.  However, the estimates in \cite{BCS} also required a
``post processing algorithm", tracing the location of the large shocks in the approximate solution.
We would like to achieve error estimates based solely on an a posteriori 
bound of the total variation.
The possibility of such estimates is the content of the following corollary.

\begin{corollary}
Let (\ref{1}) be an $n\times n$ strictly hyperbolic system, generating 
a Lipschitz semigroup of entropy-weak solutions on a domain $\D$ of functions with small total variation.
Then, given $T,R>0$, there exists a function $\ve\mapsto \varrho(\ve)$ 
with the following properties.
\begi
\item[(i)]
$\varrho$ is continuous, nondecreasing, with $\varrho(0)=0$. 
\item[(ii)] Let $t\mapsto u(t)\in \D$ be an $\ve$-approximate solution to 
(\ref{1}), satisfying {\bf (AL)-(P$_\ve$)} and supported inside the interval $[-R,R]$.
Then, calling $\bar u= u(0)$, one has
\bel{err4}
\bigl\|u(t)-S_t \bar u\bigr\|_{\L^1}~\leq~\varrho(\ve)\qquad\qquad \forall t\in [0,T].\eeq
\endi
\end{corollary}

{\bf Proof.} If the conclusion fails, there exists a sequence of $\ve_n$-approximate solutions $(u_n)_{n\geq 1}$, all supported inside $[-R,R]$, with $\ve_n\to 0$ but
\bel{err5}\sup_{t\in [0,T]} ~
\bigl\|u_n(t)-S_t u_n(0)\bigr\|_{\L^1}~\geq~\delta_0 ~>~0 \qquad\qquad \forall n\geq 1.\eeq
By compactness, taking a subsequence we achieve the $\L^1$-convergence
$u_n(t)\to u(t)$, uniformly for $t\in [0,T]$.   
Setting $\bar u(x) \doteq u(0,x)$, the limit function $u$ is thus an entropy weak
solution of (\ref{1})-(\ref{2}), distinct from the semigroup trajectory $S_t\bar u$.   This contradicts the uniqueness stated in 
Theorem~\ref{t:1}.
\endproof

We regard the function $\varrho(\cdot)$ as a {\bf universal convergence rate} 
for approximate BV solutions
to the hyperbolic system (\ref{1}).   
Having proved the existence of such a function, the major open problem is now to 
provide an asymptotic estimate on $\varrho(\ve)$,
as $\ve\to 0$.
  In some sense,
starting from a uniqueness theorem and deriving a uniform convergence rate
is a task analogous to the derivation of  quantitative compactness estimates \cite{AGN12, AGN15, AGN19, DLG}.
Based on the  convergence estimates already available 
for the Glimm scheme \cite{AM1, BM98}
and for  vanishing viscosity approximations \cite{BHWY, BY2},
one might guess that 
$\varrho(\ve)\approx \sqrt\ve |\ln \ve|$. 
We leave this as an  open question for future investigation.

{\bf Acknowledgment.} The research by the first author
 was partially supported by NSF with
grant  DMS-2006884, ``Singularities and error bounds for hyperbolic equations".
The second author acknowledges the hospitality of the
Department of Mathematics, Penn State University -- March 2023.

\end{document}